\documentclass[brochure,12pt]{bourbaki}
\usepackage[matrix,arrow]{xy}
\usepackage{amssymb,amsfonts,amsmath,footnote}
\usepackage[francais]{babel}
\usepackage[latin1]{inputenc}
\newcommand{\R}{\mathbb{R}}
\newcommand{\C}{\mathbb{C}}
\newcommand{\Q}{\mathbb{Q}}
\newcommand{\Z}{\mathbb{Z}}

\newcommand{\ra}{\rightarrow}

\date{Janvier 2013}
\bbkannee{65\`eme ann\'ee, 2012-2013}
\bbknumero{1063}
\title{Progrès récents sur les fonctions normales}
\subtitle{d'apr\`es Green-Griffiths, Brosnan-Pearlstein, M. Saito, Schnell...}
\author{Fran\c cois CHARLES}
\address{Universit\'e de Rennes 1\\
Institut de Recherche Math\'ematique de Rennes\\
UMR 6625 du CNRS\\
263 avenue du G\'en\'eral Leclerc\\
CS 74205\\
35042 Rennes Cedex}
\email{francois.charles@univ-rennes1.fr}

\begin{document}
\maketitle

\noindent{\bf INTRODUCTION}

\bigskip

Ce texte a pour but de décrire un certain nombre d'avancées récentes en théorie de Hodge concernant les valeurs limites de fonctions normales et les dégénérescences de jacobiennes intermédiaires. 

Soit $X$ une variété projective lisse sur le corps des nombres complexes. Généralisant l'application d'Abel-Jacobi habituelle associant à un diviseur homologiquement trivial sur $X$ un point de la jacobienne de $X$, Griffiths a construit dans \cite{Gr68} une application d'Abel-Jacobi pour les cycles de codimension quelconque homologues à zéro à valeurs dans un tore complexe, la jacobienne intermédiaire de $X$. Quand $X$ varie dans une famille lisse au-dessus d'une base $B$, les jacobiennes intermédiaires varient elles aussi en une famille lisse, et les familles de cycles fournissent des sections de cette fibration en jacobiennes intermédiaires. Ces sections, qui vérifient une équation différentielle venant de la propriété de transversalité de Griffiths, sont des cas particuliers de fonctions normales.

La construction qui précède est de nature transcendante : les jacobiennes intermédiaires ne sont pas en général des variétés abéliennes, et les fonctions normales sont donc des fonctions seulement analytiques. En outre, pour les problèmes de nature globale, il est nécessaire de pouvoir travailler dans un contexte où la variété $X$ ci-dessus peut dégénérer en une variété singulière. Se pose alors la question du comportement à la limite de la fibration en jacobiennes intermédiaires et des fonctions normales. 

Dans le cas d'un pinceau de Lefschetz -- donc d'une base de dimension $1$ -- ces préoccupations apparaissent dès la preuve par Poincaré de la conjecture de Hodge pour les diviseurs. En direction de la conjecture de Hodge, le théorème de Zucker sur les fonctions normales \cite{Zu76} donne une correspondance entre classes de Hodge dans la cohomologie primitive de degré $2d$ d'une variété projective lisse complexe de dimension $2d$ et certaines fonctions normales pour une compactification de la fibration en jacobiennes intermédiaires associée à un pinceau de Lefschetz. 

Pour les cycles de codimension au moins $2$, l'application d'Abel-Jacobi n'est pas surjective en général. La détermination de son image est aujourd'hui une question largement ouverte. Il est bien connu que ce problème est un obstacle majeur dans l'étude de la conjecture de Hodge par la méthode du paragraphe précédent. 

\bigskip

Depuis quelques années, à la suite d'un article de Thomas \cite{Th05} puis d'un article, fondamental pour les questions dont nous traitons, de Green et Griffiths \cite{GG07}, est apparue l'idée de construire des cycles algébriques en travaillant non pas avec un pinceau de Lefschetz, mais avec la famille universelle des sections hyperplanes d'une variété. Comme Green et Griffiths l'ont mis en avant, travailler au-dessus d'une base de grande dimension permet de faire apparaître des invariants, les singularités d'une fonction normale, qui sont de torsion dans le cas d'une base de dimension $1$. L'existence de singularités qui ne sont pas de torsion est prédite par la conjecture de Hodge, et lui est en fait équivalente.

Suivant ce cercle d'idées, l'étude des dégénérescences des fonctions normales au-dessus d'une base quelconque s'est développée de manière importante en quelques années. Un énoncé particulièrement frappant dans ce contexte est le suivant.

\begin{theo}\label{algebrique}
Soit $\nu$ une fonction normale admissible sur une variété algébrique complexe $B$. Alors le lieu d'annulation de $\nu$ est algébrique.
\end{theo}

Il s'agit d'un résultat de théorie de Hodge. Dans le théorème précédent, une fonction normale admissible est la version abstraite d'une fonction normale venant de la géométrie comme plus haut. L'algébricité du lieu des zéros ci-dessus contraste avec le fait que la définition même de fonction normale n'a de sens que dans un cadre analytique et que le lieu des zéros de $\nu$ n'est a priori qu'un sous-ensemble analytique de $X$. Nous expliquerons dans le cours du texte le lien entre les méthodes de Green-Griffiths, ainsi que les conjectures usuelles sur les cycles algébriques, et le théorème \ref{algebrique}.

D'après un théorème de Carlson \cite{Ca87}, les jacobiennes intermédiaires paramètrent des familles de structures de Hodge mixtes. Plus précisément, si $H$ est une structure de Hodge polarisée de poids  strictement négatif, la jacobienne intermédiaire $J(H)$ est en bijection canonique avec les structures de Hodge mixtes $H'$ qui sont extension de $\mathbb Z$ par~$H$. Une telle extension est scindée si et seulement si $H'$ contient une classe de Hodge -- c'est-à-dire un élément de $H'\cap W_0 H'_{\mathbb Q}\cap F^0 H_{\mathbb C}$, où $W$ est la filtration par le poids et $F$ la filtration de Hodge.  

Via cette interprétation, le théorème \ref{algebrique} est un résultat d'algébricité du lieu des classes de Hodge dans l'espace total d'une variation de structures de Hodge mixte. En ce sens, il est l'analogue d'un résultat célèbre de Cattani, Deligne et Kaplan \cite{CDK} sur l'algébricité du lieu des classes de Hodge dans l'espace total d'une variation de structures de Hodge pures polarisées de poids pair. De même que le résultat de Cattani-Deligne-Kaplan serait une conséquence de la conjecture de Hodge, le théorème \ref{algebrique} est relié aux conjectures de Bloch et Beilinson sur l'existence d'une filtration sur les groupes de Chow des variétés algébriques.

\bigskip

Pour démontrer le théorème \ref{algebrique}, il suffit, c'est le théorème GAGA de Serre \cite{Se56}, de montrer que l'adhérence du lieu des zéros de $\nu$ dans une compactification algébrique de $X$ est encore un sous-ensemble analytique. Dans une série d'articles \cite{BP1, BP2, BP3}, Brosnan et Pearlstein sont parvenus à prouver ce résultat par une analyse fine de la théorie de Schmid décrivant en un certain sens les limites de familles de structures de Hodge sur des produits de disques épointés, ainsi que de ses généralisations dues en particulier à Cattani, Deligne et Kaplan. Leur méthodes n'étudient que le lieu des zéros et ne décrivent pas la dégénérescence de la fibration en jacobiennes intermédiaires.

L'étude de cette dégénérescence, et la recherche d'un modèle de Néron pour les familles de jacobiennes intermédiaires -- en un sens à préciser plus tard -- a été entreprise par Green, Griffiths et Kerr \cite{GGK}, puis Brosnan, Pearlstein et Saito \cite{BPS}, généralisant des travaux antérieurs de Zucker \cite{Zu76} et Clemens \cite{Cle83}. Dans cette direction, le résultat le plus puissant est dû à Schnell \cite{Sc09} qui construit un modèle de Néron pour les familles de jacobiennes intermédiaires au-dessus d'une compactification lisse quelconque de la base. 

L'apport essentiel de la construction de Schnell est d'utiliser de manière efficace la théorie des modules de Hodge mixtes de Saito \cite{Sa88, Sa90}. Ces derniers, qui jouent en théorie de Hodge le rôle des faisceaux pervers $\ell$-adiques pour les variétés sur les corps finis, fournissent ici des faisceaux cohérents étendant les fibrés vectoriels qui permettent de définir la fibration en jacobiennes intermédiaires. C'est eux aussi qui permettent à la construction de fonctionner au-dessus d'une base quelconque, quand la théorie de Schmid et ses généralisations ne sont en général adaptées qu'à la situation où l'on travaille sur un produit de disques épointés.

Le modèle de Néron construit par Schnell est suffisant pour donner une démonstration du théorème \ref{algebrique}. C'est celle-là que nous esquisserons.

Il existe une troisième approche à ce type de questions, dont nous ne parlerons pas, développée principalement par Kato, Nakayama et Usui \cite{KNU08, KNU10, KNU11, KU09}. Dans cette série d'articles, les auteurs montrent comment appliquer des méthodes de géométrie logarithmique à l'étude des dégénérescences de structures de Hodge mixtes. Dans le cas où la dégénérescence se produit au-dessus d'un diviseur à croisements normaux, on peut construire un modèle de Néron qui permette de démontrer le théorème \ref{algebrique}.

\bigskip

Les résultats de dégénérescences évoqués ci-dessus font tous appel de manière fine à des calculs précis sur les structures de Hodge limites. Les techniques introduites par Schmid dans \cite{Sc73} notamment orbites nilpotentes, théorèmes de l'orbite SL$_2$, ainsi que leurs variantes et leurs généralisations, jouent un rôle dans les trois démonstrations mentionnées plus haut. Il s'agit dans tous les cas d'obtenir des estimées souvent délicates sur les normes de certaines solutions d'équations différentielles -- correspondant aux sections du système local entier ou réel sous-jacent à une variation de structures de Hodge -- au voisinage d'une dégénérescence. 

Dans ce texte, nous avons choisi d'évoquer un autre aspect du cercle d'idées en question en nous concentrant sur le versant {\og pervers\fg} de la discussion. Le formalisme des faisceaux pervers et des modules de Hodge mixtes joue un rôle crucial à la fois dans l'étude des singularités des fonctions normales et dans la construction du modèle de Néron de Schnell. 

C'est précisément cet aspect-là qui est suffisamment fonctoriel pour travailler au-dessus d'une base quelconque. En ce qui concerne les applications aux cycles algébriques, cela permet notamment de travailler avec la famille universelle des sections hyperplanes d'une variété projective lisse donnée, cas où le diviseur paramétrant le lieu des hypersurfaces singulières est loin d'être à croisements normaux. En particulier, la famille universelle des hypersurfaces de degré donné de l'espace projectif est un exemple intéressant, et les résultats de Schnell permettent de prolonger, en un sens, l'étude de la cohomologie des hypersurfaces lisses, pour lesquelles la filtration de Hodge a été décrite par Griffiths \cite{Gr69}. Une situation similaire est décrite par le théorème du support de Ngô \cite{Ngo} qui, dans un cas très particulier, décrit les faisceaux pervers intervenant dans le théorème de décomposition pour la courbe plane universelle. En un certain sens, ces résultats précisent et élargissent l'énoncé du théorème de décomposition de Beilinson-Bernstein-Deligne-Gabber \cite{BBD} dans ce contexte.

Nous renvoyons à \cite{KP10} pour une autre exposition détaillée de la théorie récente des fonctions normales.

\bigskip

Le texte est organisé comme suit. Dans la première partie, nous rappelons brièvement la théorie classique des fonctions normales, due à Griffiths. La deuxième partie est consacrée aux questions d'algébricité de lieux de Hodge et à des variations sur l'énoncé du théorème \ref{algebrique}.
La troisième partie est consacrée à l'utilisation -- largement conjecturale -- des fonctions normales dans la construction de cycles algébriques. Après avoir décrit le cas des pinceaux de Lefschetz et le théorème de Zucker, nous introduisons la notion de singularité d'une fonction normale, qui n'apparaît, à torsion près, qu'au-dessus d'une base de dimension supérieure, et nous expliquons, suivant \cite{BFNP}, comment l'existence de suffisamment de singularités pour les fonctions normales est équivalente à la conjecture de Hodge. Enfin, la dernière partie de ce texte est consacrée à la construction d'un modèle de Néron pour les familles de jacobiennes intermédiaires, suivant Schnell. Nous décrivons l'objet obtenu dans le cas de la courbe plane universelle en utilisant le théorème du support de Ngô \cite{Ngo} et donnons les grandes lignes de la preuve du théorème \ref{algebrique}.

\bigskip

Tout au long de l'exposé, le corps de base est le corps des nombres complexes. Les groupes de cohomologie singulière que nous aurons à utiliser seront toujours considérés modulo leur sous-groupe de torsion.

\bigskip

\noindent{\bf Remerciements}. Je tiens à exprimer ma reconnaissance à Christian Schnell pour ses  nombreuses réponses à mes nombreuses questions. Je remercie Gérard Laumon de m'avoir expliqué le théorème du support de Ngô dans le cas de la famille universelle des courbes planes, ainsi que Matt Kerr et Olivier Benoist.

\section{La théorie classique}

Dans cette partie, on décrit sans démonstrations la théorie classique de l'application d'Abel-Jacobi et des fonctions normales. Sous cette forme, elle est essentiellement due à Griffiths. Nous nous concentrons sur le cas du poids $-1$, mais la plupart des énoncés valent {\it mutatis mutandis} pour les structures de Hodge de poids strictement négatif.

\subsection{Jacobiennes intermédiaires et application d'Abel-Jacobi}\label{JI}

Soit $H$ une structure de Hodge entière de poids $-1$. On note $H_{\mathbb Z}$ le groupe abélien sous-jacent, $H_{\mathbb C}$ l'espace vectoriel complexe correspondant, et $F^{\bullet}$ la filtration de Hodge sur $H_{\mathbb C}$.

Le poids étant $-1$, les sous-espaces $F^0 H_{\C}$ et son conjugué complexe $\overline{F^0H_{\mathbb C}}$ sont en somme directe, et l'on a 
$$H_{\C}=F^0H_{\C}\oplus \overline{F^0H_{\C}}.$$

Il résulte de cette décomposition que l'application canonique 
$$H_{\R}\subset H_{\C}\ra H_{\C}/F^0H_{\C}$$
est un isomorphisme d'espaces vectoriels réels. Cela signifie que  le groupe abélien $H_{\Z}$ s'identifie à un réseau de l'espace vectoriel complexe $H_{\C}\ra H_{\C}/F^0H_{\C}$, et justifie la définition suivante.

\begin{defi}
La jacobienne intermédiaire de $H$ est le tore complexe 
$$J(H)=\frac{H_{\mathbb C}}{F^0H_{\mathbb C}+H_{\mathbb Z}}.$$
\end{defi}

\begin{rema}
Le point clé de la définition ci-dessus, celui qui permet de quotienter par le groupe $H_{\Z}$ et d'obtenir un espace topologique séparé, est précisément le fait que la filtration de Hodge soit opposée à sa filtration conjuguée. Cela se traduit ici par l'égalité  $F^0 H_{\C}\cap \overline{F^0 H_{\C}}=0.$

Il s'agit là d'une propriété de nature analytique, qui n'est semble-t-il pas accessible par la géométrie algébrique seule, et dont la démonstration fait appel à la théorie des opérateurs différentiels elliptiques. Les estimées nécessaires à la construction des modèles de Néron sont justement, on le verra, celles qui permettent de garantir la propriété de séparation.
\end{rema}

Si $X$ est une variété projective lisse et $k$ un entier strictement positif, la construction précédente s'applique au cas où $H$ est le groupe de cohomologie singulière $H^{2k-1}(X, \Z(k))$ muni de sa structure de Hodge canonique. On obtient ainsi la $k$-ième jacobienne intermédiaire de $X$, notée $J^k(X)$.

Si $k=1$, cette construction fournit le tore complexe sous-jacent à la variété abélienne Pic$^0(X)$, qui paramètre les fibrés en droite homologiquement triviaux sur $X$. Si $k$ est égal à la dimension de $X$, on retrouve la variété d'Albanese, qui paramètre, elle, certaines classes d'équivalence de zéros-cycles sur $X$.

Dans le cas général, la jacobienne intermédiaire est un objet de nature plus mystérieuse. En général, si $k$ est différent de $1$ et de la dimension de $X$, il ne s'agit pas d'une variété abélienne. En effet, à l'exception de ces deux cas extrêmes, une polarisation de la structure de Hodge $H^{2k-1}(X, \Z(k))$ n'induit pas, pour des raisons de signe, une polarisation du tore complexe $J^k(X)$.

\bigskip

Dans la suite de l'exposé, nous expliquerons comment généraliser la construction des jacobiennes intermédiaires lorsque la structure de Hodge $H$ varie et peut éventuellement dégénérer. Dans ce contexte, $H_{\C}$ est remplacé par un $\mathcal D$-module holonome régulier $\mathcal M$. En général, $\mathcal M$ n'est pas cohérent comme $\mathcal O$-module, mais on dispose de sous-$\mathcal O$-faisceaux cohérents de $\mathcal M$. C'est avec ceux-ci qu'il est plus facile de construire des objets géométriques prolongeant les jacobiennes intermédiaires. Pour cette raison, nous aurons besoin d'une définition différente, mais équivalente.

\begin{prop}\label{dual}
Soit $H$ une structure de Hodge de poids $-1$. Soit $\check H=H^{\vee}(1)$ la structure de Hodge duale de $H$, normalisée pour être de poids $-1$. La jacobienne intermédiaire de $H$ est canoniquement isomorphe au tore complexe 
$$\frac{(F^0\check H_{\C})^{\vee}}{H_{\Z}}.$$
\end{prop}

\noindent{\sc Preuve} --- 
La filtration de Hodge sur $\check H_{\C}$ est donnée par la formule $F^p \check H_{\C}=$\break
$\{f : H_{\C} \ra \C | f_{|F^{-p}H_{\C}}=0\}$. Le réseau entier $H_{\Z}$ s'envoie dans $(F^0 \check H_{\C})^{\vee}$ par bidualité et restriction. Cela étant, on dispose d'un accouplement canonique 
$$(F^0 \check H_{\C})\otimes H_{\C}\ra \C$$ 
qui passe au quotient pour donner un isomorphisme
$$\frac{H_{\C}}{F^0 H_{\C}}\simeq (F^0\check H_{\C})^{\vee}$$
compatible à l'inclusion de $H_{\Z}$ des deux côtés, ce qui démontre la proposition.\qed

\begin{rema}
Tout cela vaut encore si le poids est seulement strictement négatif. Dans ce cas, les jacobiennes intermédiaires ne sont cependant pas compactes en général.
\end{rema}

\bigskip

Généralisant la situation des diviseurs et des zéros-cycles, on peut définir des applications d'Abel-Jacobi à valeurs dans les jacobiennes intermédiaires. \'Etant donné un cycle $Z$ de codimension $k$ sur $X$, c'est-à-dire une combinaison formelle à coefficients entiers de sous-schémas réduits de $X$ de codimension $k$, la classe de cohomologie de $Z$ appartient au groupe $H^{2k}(X, \Z(k))$. Si $Z$ est homologue à zéro, on peut définir de manière canonique l'image $aj(Z)$ de $Z$ dans la jacobienne intermédiaire $J^k(X)$.

On peut montrer que l'application $aj$ envoie les cycles rationnellement équivalents à zéro sur l'élément neutre de la jacobienne intermédiaire. L'application d'Abel-Jacobi passe donc au quotient pour définir une application 
$$aj : CH^k(X)_{hom}\ra J^k(X)$$
définie sur le groupe de Chow des cycles homologues à zéro dans $X$. C'est un morphisme de groupes.

Pour $k=1$, on retrouve l'application d'Abel-Jacobi pour les diviseurs. Il s'agit dans ce cas d'une bijection. Si $k$ est égal à la dimension de $X$, on obtient l'application d'Albanese. Cette application est surjective. Cependant, même dans ce cas, le noyau peut être très gros comme l'a montré Mumford \cite{Mu68}. En général, le noyau comme l'image de l'application d'Abel-Jacobi sont très mal compris.

En un sens à préciser, l'application d'Abel-Jacobi est compatible à l'action des correspondances. Dans la suite de cette section, nous décrivons les résultats de Griffiths sur son comportement quand $X$ et $Z$ varient dans une famille analytique.

\subsection{Variations de structures de Hodge}\label{VHS}

Soient $B$ une variété analytique lisse et $k$ un entier. Une variation de structures de Hodge entière de poids $k$ sur $B$ est la donnée
\begin{itemize}
\item d'un système local $H_{\Z}$ de groupes abéliens libres sur $B$,
\item d'un fibré vectoriel holomorphe $\mathcal H$ sur $B$ muni d'une connection plate $\nabla$,
\item d'une filtration décroissante $F^{\bullet}$ de $\mathcal H$ par des sous-fibrés vectoriels holomorphes.
\end{itemize}

La filtration ci-dessus est appelée la filtration de Hodge. On exige en outre les propriétés et compatibilités suivantes :
\begin{itemize}
\item le système local d'espaces vectoriels complexes $H_{\C}$ obtenu à partir de $H_{\Z}$ est le système local des sections plates de $\mathcal H$, 
\item en tout point $b$ de $B$, la filtration de Hodge sur $\mathcal H$ induit, via l'identification ci-dessus, une structure de Hodge de poids $k$ sur la fibre en $b$, $H_{\Z, b}$,
\item la filtration de Hodge satisfait la propriété de transversalité de Griffiths
$$\nabla(F^p\mathcal H)\subset F^{p-1}\mathcal H\otimes \Omega^1_{B/\C}.$$
\end{itemize}

La notion de variation de structures de Hodge est la version abstraite de l'objet obtenu lorsque l'on considère la cohomologie relative d'un morphisme projectif et lisse $\pi : \mathcal X\ra B$. Dans ce cas, Deligne a prouvé \cite{De68} que le complexe $R\pi_*\Q$ dans la catégorie dérivée est somme directe de ses objets de cohomologie $R^k\pi_*\Q$. La donnée du système local $R^k\pi_*\Q$ et du faisceau de cohomologie de de Rham relative $\mathcal H^k(\mathcal X/B)$ muni de la connexion de Gauss-Manin et de la filtration de Hodge donnent alors lieu à une variation de structures de Hodge.

Comme dans le cas des structures de Hodge, on dispose d'une notion de polarisation pour les variations de structures de Hodge mixtes. Il s'agit simplement d'un accouplement compatible à la connexion de Gauss-Manin qui induit une polarisation sur les structures de Hodge au-dessus des points de $B$. Dans la suite, même si nous ne le mentionnons pas, nous ne considérerons que des variations de structures de Hodge polarisées.

\subsection{Fonctions normales}\label{NF}

Soit $H=(H_{\Z}, \mathcal H)$ une variation de structures de Hodge polarisée sur $B$, de poids $-1$. La construction des jacobiennes intermédiaires en \ref{JI} se généralise au cas relatif comme suit. Bien que la construction et la vérification de ses propriétés soient aisées, l'importance de cette construction dans le cas où la variation de structures de Hodge mixte est remplacée par un module de Hodge nous amène à donner quelques détails.

La variation de structure de Hodge duale à $H$, normalisée pour être elle aussi de poids $-1$, a pour fibré vectoriel sous-jacent le dual $\check{\mathcal H}$ de $\mathcal H$. Comme dans la proposition \ref{dual}, on dispose d'une flèche naturelle
$$H_{\Z}\hookrightarrow (F^0\check{\mathcal H})^{\vee}$$
du système local $H_{\Z}$ dans le fibré vectoriel $(F^0\check{\mathcal H})^{\vee}$.

Soit $p_{\Z} : T_{\Z}\ra B$ l'espace étalé du faisceau $H_{\Z}$. C'est un espace analytique en groupes sur $B$ dont les sections s'identifient aux sections de $H_{\Z}$. De même, soit $T(F^0\check{\mathcal H})$ le spectre relatif de l'algèbre symétrique sur le faisceau localement libre $F^0\check{\mathcal H}$. L'espace $T(F^0\check{\mathcal H})$ est l'espace total du fibré vectoriel $F^0\check{\mathcal H}$.

\`A l'inclusion $H_{\Z} \hookrightarrow (F^0\check{\mathcal H})^{\vee}$ correspond une flèche 
$$\epsilon : T_{\Z}\ra T(F^0\check{\mathcal H})$$
au-dessus de $B$, compatible aux structures de groupes. 

\begin{prop}\label{epsilon}
Le morphisme $\epsilon$ est une immersion fermée.
\end{prop}

\noindent{\sc Preuve} --- L'injectivité de $\epsilon$ se démontre fibre par fibre, et vient de la construction des jacobiennes intermédiaires de \ref{JI}. Pour montrer que $i$ est une immersion fermée, on peut raisonner localement sur $B$ et donc supposer que $H_{\Z}$ est le système local trivial. 

Dans ce cas, $T_{\Z}$ est simplement une union dénombrable de copies de $B$ indexée par le groupe discret des sections globales de $H_{\Z}$. En particulier, la restriction de $p_{\Z}$ \`a chaque composante connexe de $T_{\Z}$ est propre, ce qui implique que la restriction de $\epsilon$ à chaque composante connexe de $T_{\Z}$ est propre. La restriction de $\epsilon$ à chacune des composantes connexes de $T_{\Z}$ est donc une immersion fermée. L'injectivité globale de $\epsilon$ permet de conclure.\qed

Ce qui précède permet de construire le quotient d'espaces analytiques en groupes $J=T(F^0\check{\mathcal H})/T_{\Z}$. Ses fibres sont les jacobiennes intermédiaires des fibres de la variation de structures de Hodge $H$. Puisque $\epsilon$ est une immersion fermée, il est facile de vérifier que $J$ est lisse sur $B$. 

\begin{defi}
L'espace analytique $J\ra B$ est la fibration en jacobiennes intermédiaires associée à $H$.
\end{defi}

Ce qui précède permet d'introduire la notion de fonction normale. Supposons d'abord donné un morphisme projectif lisse $\mathcal X\ra B$. La cohomologie de degré $2k-1$ des fibres de $\mathcal X$ fournit, après un twist de Tate convenable, la fibration en jacobiennes intermédiaires $J^k(\mathcal X)\ra B$. 

Si maintenant $\mathcal Z$ est une famille de cycles de codimension $k$ de $\mathcal X$ au-dessus de $B$, dont la restriction à chaque fibre est homologiquement triviale, l'application d'Abel-Jacobi donnée fibre par fibre comme en \ref{JI} fournit une section holomorphe $\nu$ de $J^k(\mathcal X)\ra B$. La section $\nu$ est le prototype d'une {\it fonction normale}. 

En général, les fonctions normales sont des sections holomorphes de la fibration en jacobiennes intermédiaires $J\ra B$. Comme dans le cas des variations de structures de Hodge, les fonctions normales venant de la géométrie vérifient une équation différentielle, c'est la condition d'horizontalité -- dont nous verrons plus tard qu'il s'agit d'une forme mixte de la transversalité de Griffiths.

Nous ne définissons pas la condition d'horizontalité dans ce paragraphe, préférant la repousser à la section \ref{MHS}. Les fonctions normales sont par définition les sections holomorphes horizontales de $J\ra B$.

Dans le cas, qui nous intéresse tout particulièrement, où $B$ est le complémentaire dans une variété analytique lisse d'un fermé analytique, El-Zein et Zucker \cite{EZZ84}, Kashiwara \cite{Ka86} et Saito \cite{Sa96} dans un contexte plus général, ont isolé une condition d'admissibilité pour les fonctions normales qui impose des conditions de croissance modérée à l'infini. Bien qu'elle soit importante techniquement, nous avons choisi de ne pas définir cette notion et nous parlerons librement de fonctions normales admissibles.

\section{Lieux de Hodge -- questions d'algébricité}

Dans cette section, nous essayons de motiver l'énoncé du théorème \ref{algebrique} et en donnons quelques variantes.

\subsection{Le théorème de Cattani-Deligne-Kaplan}\label{CDK}

Une référence particulièrement agréable pour ce qui vient est l'article \cite{Vo}.

Soit $B$ une variété quasi-projective lisse, et soit $\pi : \mathcal X\ra B$ un morphisme projectif et lisse. Soit $k$ un entier. Nous considérons dans ce paragraphe la variation de structures de Hodge de poids pair $2k$ associée à $R^{2k}\pi_*\Z$.

Soit $b$ un point de $B$, et soit $h$ un élément de $H^{2k}(\mathcal X_b, \Z)$. On dit que $h$ est une classe de Hodge si l'image de $h$ dans $H^{2k}(\mathcal X_b, \C)$ appartient à $F^kH^{2k}(\mathcal X_b, \C)$. Le transport parallèle associe à $h$ un élément canonique du groupe $H^{2k}(\mathcal X_{b'}, \Z)$ pour tout point $b'$ dans un petit voisinage -- simplement connexe -- de $b$ dans $B$. Le lieu de Hodge $T$ de $h$ dans $B$ est le lieu des points $b'$ comme ci-dessus où $h$ reste une classe de Hodge -- de manière équivalente, où $h$ reste dans $F^kH^{2k}$.

Il résulte du fait que la filtration de Hodge varie de manière holomorphe que $T$ est le germe d'un sous-ensemble analytique de $B$. Cependant, la conjecture de Hodge prévoit plus. Rappelons que celle-ci prédit que les classes de Hodge sont exactement, à un multiple entier près, les classes de cohomologie des cycles algébriques. Supposons-la vérifiée. Dans ce cas, le lieu de Hodge $T$ est le lieu des points $b'$ de $B$ où le transport parallèle de $h$ est une classe de cycle. Un argument simple, utilisant le théorème de Baire, montre alors qu'il existe une famille de cycles sur $\mathcal X$ au-dessus de $B$ dont la projection, au voisinage de $b$, contient $T$ comme composante irréductible. En particulier, $T$ est le germe d'une sous-variété algébrique.

C'est ce résultat que prouvent de manière inconditionnelle Cattani, Deligne et Kaplan dans \cite{CDK}.

\begin{theo}\label{pure}
Soit $H$ une variation de structures de Hodge polarisable de poids pair sur une variété algébrique complexe $B$. Les lieux de Hodge associés à $H$ sont des germes de sous-ensembles algébriques de $B$.
\end{theo}

La signification géométrique et arithmétique de l'énoncé précédent a été explorée par Voisin dans \cite{Vo07, Vo}. 

Dans ce qui suit, nous expliquons pourquoi le théorème \ref{algebrique} est l'analogue dans le cas mixte du théorème précédent, et nous expliquons pourquoi les conjectures de Bloch et Beilinson motivent le résultat d'algébricité dans ce cadre.

\subsection{Fonctions normales et variations de structures de Hodge mixtes}\label{MHS}

Nous commençons par un résultat de Carlson \cite{Ca87}. 

\begin{prop}\label{ext}
Soit $H$ une structure de Hodge entière de poids $-1$. Il existe une bijection canonique entre la jacobienne intermédiaire $J(H)$ et le groupe des extensions de structures de Hodge mixtes Ext$^1_{MHS}(H, \Z)$ qui paramètre les extensions
$$0\ra H\ra H'\ra \Z\ra 0.$$
\end{prop}

\noindent{\sc Preuve} --- Nous montrons simplement comment associer à toute extension comme ci-dessus un élément de $J(H)$. Considérons une telle suite exacte. Dualisant, et tordant par $\Z(1)$, on obtient une extension de structures de Hodge mixtes
$$0\ra \Z(1)\ra \check H' \ra \check H\ra 0.$$
Les morphismes de structures de Hodge étant stricts, on a $F^0\check H'_{\C}\simeq F^0\check H$.

Soit $h'$ un élément de $H'_{\Z}$ s'envoyant sur $1\in\Z$. L'élément $h'$ est bien défini modulo le groupe $H_{\Z}$. Il définit par bidualité un élément de  $(F^0\check H'_{\C})^{\vee}\simeq (F^0\check H)^{\vee}$, qui est donc bien défini modulo $H_{\Z}$. La proposition \ref{dual} fournit bien un élément de $J(H)$.\qed

\begin{rema}\label{extgen}
Dans ce qui précède, on pourrait remplacer $H$ par une structure de Hodge mixte entière dont tous les poids sont strictement négatifs. 
\end{rema}

La proposition précédente se combine à \ref{NF} et montre que, dans le cas où $H$ est une famille de structures de Hodge de poids $-1$ au-dessus d'une base $B$, les sections holomorphes de la vibration en jacobiennes intermédiaires $J(H)$ correspondent aux familles de structures de Hodge mixtes qui sont extension de la structure de Hodge constante $\Z$ par $H$. Ici, par famille de structures de Hodge mixtes nous entendons une donnée comme au début de \ref{VHS}, munie en outre d'une filtration croissante par le poids, et soumise aux mêmes compatibilités.

Cela permet de définir simplement la condition d'horizontalité de \ref{NF} pour une section holomorphe $\nu$ de $J(H)$. Une telle section $\nu$ définit une famille de structures de Hodge mixtes $H'$ au-dessus de $B$. On dit que $\nu$ est horizontale si $H'$ est une variation de structures de Hodge mixtes, c'est-à-dire si $H'$ vérifie la condition de transversalité de Griffiths.

Dans \cite{Ka86, Sa90}, Kashiwara et Saito -- dans le contexte des modules de Hodge mixtes -- ont dégagé une condition d'admissibilité pour les variations de structures de Hodge mixtes. Dans \cite{Sa96}, Saito montre qu'une fonction normale est admissible si et seulement si la variation de structures de Hodge mixtes correspondante l'est.

\begin{rema}\label{NFgen}
Suivant la remarque \ref{extgen}, on peut définir des fonctions normales associées à n'importe quelle variation de structures de Hodge mixtes dont les poids sont strictement négatifs. C'est dans cette généralité que le théorème \ref{algebrique} est vrai. Le cas pur est l'objet de ce texte, et nous nous concentrons sur le cas du poids $-1$. Le cas mixte est une conséquence du théorème \ref{gen} ci-dessous.
\end{rema}

\subsection{Algébricité du lieu de Hodge pour une variation admissible de structures de Hodge mixtes}

Le paragraphe précédent montre en quel sens le théorème d'algébricité \ref{algebrique} est un analogue mixte du théorème de Cattani, Deligne et Kaplan. 

\begin{defi}
Soit $H$ une structure de Hodge mixte. On dit qu'un élément $h$ de $H_{\Z}$ est une classe de Hodge si $h$ appartient à $W_0 H_{\Q}\cap F^0H_{\C}$, où $W_{\bullet}$ est la filtration par le poids et $F^{\bullet}$ est la filtration de Hodge.
\end{defi}

\begin{rema}
\`A un twist de Tate près, il s'agit des mêmes classes de Hodge qu'en \ref{CDK}.
\end{rema}

\begin{rema}
Par définition, les classes de Hodge dans $H$ s'identifient aux morphismes de structures de Hodge mixtes de la structure de Hodge $\Z$ dans $H$.
\end{rema}

De manière tautologique, on trouve :

\begin{prop}
Soit $H$ une structure de Hodge de poids $-1$. Soit $x$ un point de $J(H)$ et soit
$$0\ra H\ra H'\ra \Z\ra 0$$
l'extension de structures de Hodge mixtes associée. Alors $x$ est nul dans $J(H)$ si et seulement si $H'$ contient une classe de Hodge qui s'envoie sur $1$ dans $\Z$. De même, $x$~est un point de torsion de $J(H)$ si et seulement si $H'$ contient une classe de Hodge non nulle.
\end{prop}

Ce qui précède montre que le lieu des zéros d'une fonction normale est précisément une union dénombrable de lieux de Hodge pour la variation de structures de Hodge mixtes associée. Le théorème \ref{algebrique} est en ce sens une version du théorème \ref{pure} pour certaines variations de structures de Hodge mixtes admissibles qui apparaissent par extension de $\Z$ par une variation pure de poids $-1$. En fait, la combinaison de \ref{algebrique} et \ref{pure} donne en toute généralité le résultat suivant. 

\begin{theo}\label{gen}
Soit $H$ une variation de structures de Hodge mixtes sur une variété algébrique complexe $B$. Supposons $H$ admissible et polarisable. Alors le lieu de Hodge pour $H$ est une réunion dénombrable de sous-ensembles algébriques de $B$.
\end{theo}

\noindent{\sc Preuve} --- Sans perte de généralité, on peut supposer que la graduation par le poids est définie sur $\Z$. Quitte à remplacer $H$ par $W_0 H$, on se ramène au cas où les poids de $H$ sont tous négatifs.

Soient maintenant $b$ un point de $B$ et $h$ une classe de Hodge dans $H_b$. L'image de $h$ dans Gr$^W_0 H_b=W_0 H/W_{-1} H_b$ est une classe de Hodge dans la structure de Hodge pure Gr$^W_0 H_b$. Par \ref{pure}, le lieu de Hodge de l'image de $h$ dans Gr$^W_0 H_b$ est un sous-ensemble algébrique de $B$. Quitte à restreindre à ce dernier, on peut donc supposer que l'image de $h$ Gr$^W_0 H_b$ est une classe de Hodge au-dessus de $B$ tout entier. 

L'action de monodromie sur Gr$^W_0 H_b$ respecte une polarisation, qui est définie (positive ou négative) sur l'espace des classes de Hodge. Cela implique que l'orbite de $h$ dans Gr$^W_0 H_b$ sous la monodromie est finie. On peut donc supposer que $h$ est invariant par la monodromie, et correspond donc à un morphisme de la famille de structures de Hodge constante $\Z$ dans Gr$^W_0 H$.

On dispose d'une extension de structures de Hodge mixtes
$$0\ra W_{-1} H \ra H \ra \mathrm{Gr}^W_0 H\ra 0.$$
Tirant en arrière par le morphisme $\Z\ra \mathrm{Gr}^W_0 H$ induit par $h$, on peut maintenant supposer que $H$ se place dans une suite exacte 
$$0\ra W_{-1} H \ra H \ra \Z\ra 0.$$

Par dévissage et par récurrence sur les poids qui apparaissent dans $W_{-1} H$, on se ramène facilement au cas où $W_{-1}  H$ est une variation admissible de structures de Hodge qui est pure (de poids strictement négatif). Le théorème \ref{algebrique} permet de conclure grâce à \ref{MHS} et la remarque \ref{NFgen}. \qed

\subsection{Remarques sur l'algébricité et le noyau de l'application d'Abel-Jacobi}

Dans ce paragraphe, on se place dans un cadre géométrique. Nous discutons brièvement de quelques liens, mis d'abord en évidence par Bloch et Beilinson, entre certaines conjectures générales sur les groupes de Chow et le théorème \ref{algebrique}. Nous renvoyons aux articles \cite{Be87, Bl10, Ja94} pour des détails sur ces conjectures. Certains aspects de ce lien sont abordés dans \cite{Ch10}.

Soit d'abord $X$ une variété projective lisse sur $\C$. Bloch et Beilinson ont conjecturé l'existence d'une filtration décroissante $F^{\bullet}$ sur les groupes de Chow $CH^k(X)_{\Q}$ de $X$ à coefficients dans $\Q$. Cette filtration permet d'interpréter le théorème de Mumford \cite{Mu68} sur le noyau de l'application d'Abel-Jacobi pour les zéros-cycles sur les surfaces. Rappelons que ce dernier montre que l'existence d'une structure de Hodge non-triviale sur le second groupe de cohomologie d'une surface implique que ce noyau est trop gros pour être param\'etré par une variété algébrique.

Une filtration de Bloch-Beilinson sur $CH^k(X)_{\Q}$ doit être finie, fonctorielle, compatible aux correspondances, et doit vérifier notamment que le gradué $Gr_F^i Ch^k(X)_{\Q}$ doit être contrôlé, en un certain sens, par le groupe de cohomologie $H^{2k-i}(X, \Q)$. Nous renvoyons à \cite{Ja94} pour une discussion détaillée des propriétés d'une telle filtration.

La filtration de Bloch-Beilinson devrait provenir de la suite spectrale
$$E_2^{p,q}=\mathrm{Ext}^p_{MM(\C)}(\mathbf 1, \mathfrak h^q(X)(k))\implies \mathrm{Hom}_{D^b(MM(\C))}(\mathbf 1, \mathfrak h(X)(k)[p+q]).$$
Dans ce qui précède, $MM(\C)$ est la catégorie (conjecturale) des motifs mixtes sur $\C$. Pour $p+q=2i$, l'aboutissement de cette suite spectrale est le groupe de Chow $CH^k(X)$ par un théorème de Voevodsky \cite{MVW06}. Pour des raisons de poids, la suite spectrale devrait dégénérer en $E_2\otimes \Q$, d'où la filtration de Bloch-Beilinson. 

Si le terme $F^1 CH^k(X)_{\Q}$ est clair -- c'est le noyau de l'application classe de cycles à valeurs dans $H^{2k}(X, \Q(k))$ -- le terme $F^2 CH^k(X)_{\Q}$ est moins bien compris. Il semble raisonnable de s'attendre à ce que ce terme soit exactement le noyau de l'application d'Abel-Jacobi. C'est en tout cas ce que suggère le fait que $Gr_F^1 Ch^k(X)_{\Q}=F^1 CH^k(X)_{\Q}/F^2 CH^k(X)_{\Q}$ ne dépende que de la cohomologie de $X$ en degré $2k-1$.

Si c'est le cas, et puisque la filtration de Bloch-Beilinson est unique si elle existe, ce noyau doit être invariant sous l'action des automorphismes de $\C$. Cela étant, un argument élémentaire montre que, dans la situation géométrique, le lieu des zéros d'une fonction normale doit être un sous-ensemble algébrique de la base. Cela se fonde sur le fait qu'une réunion dénombrable de sous-ensembles analytiques d'une variété algébrique globalement invariante par les automorphismes du corps $\C$ est en fait réunion de sous-ensembles algébriques.

En fait, cet argument ne repose que sur la propriété suivante. On a vu plus haut que l'annulation d'une fonction normale était causée par l'existence de certaines classes de Hodge dans une variation de structures de Hodge mixtes. Dans le cas géométrique, cette dernière vient elle aussi de la géométrie, et vient donc avec une version étale $\ell$-adique. On peut dans ce cas se demander si les classes de Hodge en question sont absolues suivant la définition de Deligne \cite{De82}, voir aussi \cite{CS11}. Si c'est le cas, l'argument précédent s'applique et montre l'algébricité du lieu des zéros.

\bigskip

Pour conclure cette partie, on peut remarquer que le théorème \ref{algebrique}, ainsi que la version du théorème \ref{pure} de Cattani, Deligne et Kaplan, prouvent en fait un  résultat plus fort que ce qui est prévu par les conjectures usuelles sur les cycles algébriques. En effet, le théorème \ref{algebrique} montre non seulement que les composantes du lieu des zéros d'une fonction normale sont algébriques, mais il montre en outre que ces composantes sont en nombre fini. De même, le théorème \ref{gen} admet une version plus forte dans laquelle est prouvée une condition de finitude en fonction d'une polarisation.

Comme expliqué dans \cite[7.3]{Vo}, les contre-exemples à la conjecture de Hodge entière rendent cet énoncé de finitude pour le théorème \ref{pure} assez mystérieux. De la même façon, nous ne connaissons pas d'explication motivique de la finitude du nombre de composantes du lieu des zéros d'une fonction normale venant de la géométrie.

\section{Fonctions normales et cycles algébriques}

\subsection{Le théorème des fonctions normales et la conjecture de Hodge pour les diviseurs}\label{ZP}

Dans ce paragraphe, nous décrivons la méthode qu'utilise Poincaré pour démontrer la conjecture de Hodge pour les diviseurs, ainsi que sa généralisation suivant Zucker \cite{Zu76}. Une référence pour les classes de cohomologie de fonctions normales est \cite{Vo02}.

Soient $\overline B$ et $\overline{\mathcal X}$ deux variétés quasi-projectives lisses, et soit $\pi : \mathcal X\ra \overline B$ un morphisme plat. Soient $B$ l'ouvert de $\overline B$ au-dessus duquel $\pi$ est lisse, et $\mathcal X$ l'ouvert de $\overline{\mathcal X}$ correspondant. Dans ce paragraphe, on se restreint au cas où la base $\overline B$ est de dimension $1$. On suppose la dimension de $\mathcal X$ paire et égale à $2n$.

Soit $k$ un entier, et soit $\alpha\in H^{2n}(\overline{\mathcal X}, \Z(n))$ une classe de Hodge. On suppose que la restriction de $\alpha$ aux fibres lisses de $\pi$ est triviale. La fibration $J^n(\mathcal X/B)\ra B$ admet une sous-fibration constante venant de la jacobienne intermédiaire de l'espace total. Notons $J\ra B$ la fibration quotient. Il s'agit de la fibration en jacobiennes intermédiaires associée au quotient de la cohomologie des fibres en degré $2n-1$ par l'image de la cohomologie de l'espace total $\overline{\mathcal X}$.

On peut associer à la classe de cohomologie $\alpha$ une fonction normale qui est une section de $J$. Cela peut se voir de la manière suivante. Au-dessus de $B$, on dispose de la variété $\mathcal X\times B$ munie de la seconde projection. Le morphisme 
$$\mathrm{Id}\times \pi : \mathcal X\ra \mathcal X\times B$$
est une immersion fermée au-dessus de $B$. Soit $\psi : \mathcal U\ra B$ l'ouvert complémentaire. Via la suite exacte ouvert-fermé complémentaire relative et le fait que la partie de poids $k$ de la cohomologie de $\mathcal X$ en degré $k$ est l'image de la cohomologie de $\overline X$, on obtient une suite exacte de structure de Hodge mixtes
$$H^{2n-1}(\overline{\mathcal X}, \Z(n))\otimes\Z_B \ra R^{2n-1}\pi_* \Z(n)\ra R^{2n-1}\psi_* \Z(n)\ra H^{2n}(\overline{\mathcal X}, \Z(n))\otimes\Z_B,$$
où $\Z_S$ est la structure de Hodge constante sur $B$. Tirant en arrière par la classe de $\alpha$, qui est dans l'image de la dernière flèche car elle se restreint à zéro sur les fibres de $\pi$, on obtient une extension de structures de Hodge mixtes qui définit, par \ref{MHS}, une fonction normale $\nu_{\alpha}$ pour $J$. Il est montré dans \cite{Sa96} que $\nu_{\alpha}$ est admissible.

\bigskip

Supposons maintenant que $\pi$ soit le morphisme obtenu à partir d'un pinceau de Lefschetz de sections hyperplanes sur une variété projective lisse $X$ en éclatant le lieu de base du pinceau. Pour fixer les idées et simplifier les notations, supposons en outre que la cohomologie de $\overline{\mathcal X}$ soit nulle en degré $2n-1$. Dans ce cas, $J$ est la fibration $J^n(\mathcal X)$ et $\nu_{\alpha}$ correspond à une extension de structures de Hodge mixtes 
$$0\ra R^{2n-1}\pi_* \Z(n)\ra H'\ra \Z\ra 0.$$

\begin{defi}\label{classe}
La classe de cohomologie $[\nu_{\alpha}]$ de $\nu_{\alpha}$ est l'image de $1$ dans le groupe $H^1(B, R^{2n-1}\pi_*\Z(n))$ par la suite exacte longue venant de l'extension précédente.
\end{defi}

\begin{rema}
Cette définition vaut telle quelle pour une fonction normale associée à une variation de structure de Hodge de poids $-1$ quelconque.
\end{rema}

On peut montrer, en chassant les diagrammes, que la classe $[\nu_{\alpha}]$ s'obtient à partir de $\alpha$ par la suite spectrale de Leray. Grâce au théorème de Lefschetz fort, cela montre que $\nu_{\alpha}$ détermine $\alpha$. En outre, si $\alpha$ est la classe de cohomologie d'un cycle $Z$, on peut montrer que $\nu_{\alpha}$ est bien la fonction normale associée à $Z$ par l'application d'Abel-Jacobi relative.

\bigskip

Ce qui précède suggère la méthode suivante pour démontrer la conjecture de Hodge pour $\alpha$ : on peut espérer, dans certains cas, que les jacobiennes intermédiaires des fibres de $\pi$ paramètrent des cycles algébriques sur les fibres. La fonction $\nu_{\alpha}$ paramètre alors une famille (a priori holomorphe) de cycles sur les fibres de $\pi$. L'espace total de cette famille est alors un candidat pour un cycle dont la classe de cohomologie serait $\alpha$.

Cette stratégie se heurte à deux problèmes. Le premier est d'assurer effectivement que les jacobiennes intermédiaires paramètrent des cycles algébriques. Cela est faux en général, dès que $n$ est strictement supérieur à $1$. C'est un obstacle majeur pour la conjecture de Hodge, et correspond au fait que l'image de l'application d'Abel-Jacobi pour les cycles de codimension supérieure est très mal connue. 

Même dans le cas où le problème ci-dessus ne se pose pas, en particulier dans le cas de la conjecture de Hodge pour les diviseurs sur les surfaces, le problème se pose de montrer que la famille holomorphe de cycles ainsi obtenues est en fait algébrique. Par \cite{Se56}, ce problème est essentiellement équivalent à celui de prolonger la fibration $J$ et la section $\nu_{\alpha}$. Sans hypothèse sur $n$, c'est ce qui est réalisé dans \cite{Zu76, EZZ84} -- toujours dans le cas d'un pinceau de Lefschetz.

Dans cette situation, on peut construire une compactification de $J$ au-dessus de $\overline B$ de la manière suivante. Soit $\overline{\mathcal H}$ l'extension canonique de Deligne \cite{De70} du fibré vectoriel $\mathcal H = \mathcal H^{2n-1}(\mathcal X/B)(n)$ de cohomologie de de Rham relative. Il s'agit d'un fibré vectoriel canonique tel que la connexion de Gauss-Manin sur $\mathcal H$ s'étende en une connexion à pôles logarithmiques sur $\overline{\mathcal H}$.

La filtration de Hodge s'étend à $\overline{\mathcal H}$ -- c'est élémentaire dans ce cas car $\overline B$ est de dimension $1$. Soit $T_{\Z}$ l'espace étalé du faisceau $j_*R^{2n-1}\pi_* \Z(n)$, où $j$ est l'inclusion de $B$ dans $\overline B$, et soit $T$ l'espace total du fibré $\overline{\mathcal H}/F^0\overline{\mathcal H}$. On peut montrer que l'on a une injection $T_{\Z}\hookrightarrow T$ au-dessus  de $\overline B$ prolongeant l'inclusion de la proposition \ref{epsilon}.

\begin{defi}
L'extension de Zucker de la famille de jacobiennes intermédiaires $J\ra B$ est le quotient
$$J^Z(R^{2n-1}\pi_* \Z(n))=T/T_{\Z}.$$
\end{defi}

Les deux résultats qui suivent correspondent au \og théorème sur les fonctions normales\fg\,de Zucker.

\begin{prop}\label{sepZ}
L'espace analytique $J^Z(R^{2n-1}\pi_* \Z(n))$ est séparé.
\end{prop}

\begin{rema}
La définition ci-dessus est valable en toute généralité, au-dessus notamment d'une base de dimension quelconque, pourvu que le complémentaire de $B$ dans $\overline B$ soit un diviseur à croisements normaux. Cependant, la séparation de l'extension de Zucker, cruciale pour les applications, n'est vraie que dans le cas d'un pinceau de Lefschetz. En effet, elle est fausse en général \cite{Sa96, GGK} mais vaut dès que les monodromies locales $T$ satisfont $(T-1)^2=0$ comme prouvé dans \cite{Zu76}.
\end{rema}

\begin{prop}\label{extZ}
Il existe un entier non nul $m$ tel que la fonction normale $m\nu_{\alpha}$ se prolonge en une section de l'extension de Zucker de $J$.
\end{prop}

Nous préciserons la signification géométrique de l'entier $m$ un peu plus bas.

Ces deux résultats permettent de démontrer la conjecture de Hodge pour les diviseurs par la méthode décrite ci-dessus : pour $n=1$, on obtient, après multiplication par $m$, une famille de diviseurs de degré zéro sur les fibres de $\pi$ qui s'étend de manière holomorphe -- c'est la conjonction de la propriété d'extension de la fonction normale et de la séparation de $J^Z$. Par \cite{Se56}, la famille est algébrique, et son espace total fournit un diviseur de $\mathcal X$, qui est celui prédit par la conjecture de Hodge via les diverses compatibilités mentionnées ci-dessus.

\subsection{Singularités des fonctions normales}

Nous introduisons ici la notion de singularité d'une fonction normale, suivant \cite{GG07, dCM09, BFNP}. Il s'agit d'un invariant qui n'apparaît qu'au-dessus d'une base de dimension au moins $2$, pour des raisons que nous expliquerons.

Soient $\overline B$ une variété analytique lisse, $B$ un ouvert de $\overline B$, $H$ une variation de structures de Hodge de poids $-1$ et $\nu$ une fonction normale admissible pour $J=J(H)$. On note comme plus haut $H_{\Z}$ le système local sous-jacent à $H$.

Soit $d$ la dimension de $B$, et notons $j$ l'inclusion de $B$ dans $\overline B$. Si $E$ est un système local sur $B$ et si $b$ est un point de $\overline B$, on note $H^k(E)_b$ le groupe limite projective des $H^k(U\cap B, E)$, où $U$ parcourt les voisinages ouverts de $b$ dans $\overline B$. Autrement dit, on a 
$$H^k(E)_b=H^k(\{b\}, i^*Rj_*E)$$
où $i$ est l'inclusion de $\{b\}$ dans $\overline B$.

La définition suivante est due à Green et Griffiths dans \cite{GG07}.
\begin{defi}
La singularité $\sigma_{\Z,b}(\nu)$ de $\nu$ en $b$ est l'image de $[\nu]$ dans $H^1(H_{\Z})_b$. On note $\sigma_b(\nu)$ son image dans $H^1(H_{\Q})_b$.
\end{defi}

Ici $[\nu]$ est la classe de cohomologie de $\nu$ définie en \ref{classe}. La singularité de $\nu$ mesure la non-trivialité topologique au voisinage de $b$ de l'extension de systèmes locaux sous-jacente à l'extension de structures de Hodge mixtes correspondant à $\nu$. Bien entendu, elle est nulle si $b$ appartient à $B$.

\bigskip

C'est dans le contexte des faisceaux pervers qu'il faut envisager la notion de singularité d'une fonction normale. Cela nécessite en particulier de travailler avec des coefficients rationnels. Nous ne rappelons pas les éléments de la théorie et renvoyons à \cite{BBD, dCM} pour la notion de faisceau pervers et pour le théorème de décomposition de Beilinson-Bernstein-Deligne-Gabber. 

Le système local $H_{\Q}$ définit un faisceau pervers $H_{\Q}[d]$. Son extension intermédiaire $j_{!*}H_{\Q}[d]$ est un faisceau pervers, c'est le complexe d'intersection $IC(H_{\Q})$. On note $IH^k(H_{\Q})_b$ la cohomologie d'intersection associée en $b$. On a par définition 
$$IH^k(H_{\Q})_b=H^{k-d}(\{b\}, i^*IC(H_{\Q})),$$
d'où une flèche
$$IH^k(H_{\Q})_b\ra H^k(H_{\Q})_b.$$
On peut montrer que cette flèche est injective pour $k=1$.

Notre but est maintenant d'esquisser la preuve du résultat suivant, démontré dans \cite{BFNP, dCM09}. 

\begin{theo}\label{perv}
La singularité $\sigma_b(\nu)$ appartient à l'image de $IH^1(H_{\Q})_b$ dans $H^1(H_{\Q})_b$.
\end{theo}

On trouve en particulier :

\begin{coro}\label{torsion}
Le lieu des points $b$ de $\overline B$ tel que $\sigma_b(\nu)$ est non nul est de codimension au moins $2$ dans $\overline B$. En particulier, si $B$ est de dimension $1$, la singularité $\sigma_{\Z, b}(\nu)$ est de torsion pour tout point $b$ de $\overline B$.
\end{coro}

\noindent{\sc Preuve du corollaire} ---
Puisque $IC(H_{\Q})$ est un faisceau pervers, le support de $\mathcal H^{1-d}(IC(H_{\Q}))$ est de codimension au moins $2$. On conclut par l'injectivité de $IH^1(H_{\Q})_b\ra H^1(H_{\Q})_b.$\qed

\begin{rema}
La notion de singularité permet de comprendre les questions de prolongement de fonctions normales à l'extension de Zucker comme en \ref{ZP} : une fonction normale (admissible) se prolonge au modèle de Zucker si et seulement si ses singularités sont nulles. En particulier, l'entier $m$ de la proposition \ref{extZ} est le pgcd des ordres des singularités locales. Cela vaut encore dans le cas où le complémentaire de $B$ dans $\overline B$ est un diviseur à croisements normaux. Dans ce cas, la singularité peut n'être pas de torsion et aucun multiple de la fonction normale ne se prolonge à l'extension de Zucker.
\end{rema}

L'existence de singularités qui ne sont pas de torsion est prédite par la conjecture de Hodge. Nous le verrons plus bas.

\bigskip

\noindent{\sc Preuve du théorème} ---
Nous ne donnons que les grandes lignes de la démonstration. Soit 
$$0\ra H_{\Q}\ra H'\ra \Q\ra 0$$
l'extension de systèmes locaux de $\Q$-espaces vectoriels sur $B$ induite par la fonction normale $\nu$.

Considérons le complexe de faisceaux pervers
$$0\ra IC(H_{\Q})\ra j_{!*}H'[d]\ra\Q[d]\ra 0$$
obtenu par extension intermédiaire. On va montrer que ce complexe est exact. Pour des raisons générales, le foncteur $j_{!*}$ préserve injections et surjections de faisceaux pervers. C'est une conséquence formelle de l'exactitude à gauche de $^pj_*$ et de l'exactitude à droite de $^pj_!$.

Pour montrer l'exactitude de tout le complexe ci-dessus, il faut raisonner avec des poids. Comme on le verra dans la section suivante, c'est le formalisme des modules de Hodge mixtes de Saito qui est adapté dans ce cadre. Dans le cas géométrique, où $\nu$ est une fonction normale venant d'une famille de cycles homologues à zéro sur une base algébrique, \cite{BBD} fournit, après tensorisation par $\Q_{\ell}$ et par spécialisation de la situation à un corps fini, la notion de poids nécessaire. Nous parlerons donc librement de poids et de pureté pour $H_{\Q}$ et $H'$, au prix d'un abus de langage.

Le système local $H_{\Q}[d]$ est pur de poids $-1$ par hypothèse, et $\Q[d]$ est pur de poids $0$. Par conséquent, le complexe d'intersection $IC(H_{\Q})$ est pur de poids $-1$. Les poids de $j_{!*}H'[d]$ sont donc $0$ et $-1$. Considérons le gradué Gr$^W_0 j_{!*}H'[d]$. La flèche Gr$^W_0 j_{!*}H'[d]\ra \Q[d]$ est surjective comme mentionné plus haut, ce qui permet d'écrire 
$$\mathrm{Gr}^W_0 j_{!*}H'[d]=\Q[d]\oplus D$$
où $D$ est supporté sur le complémentaire de $B$ dans $\overline B$. Mais on dispose d'une surjection 
$$ j_{!*}H'[d]\ra \mathrm{Gr}^W_0 j_{!*}H'[d]\ra D.$$
Comme $j_{!*}H'[d]$ n'a pas de quotient non trivial supporté sur le complémentaire de $B$ dans $\overline B$, on trouve $D=0$ et $\mathrm{Gr}^W_0 j_{!*}H'[d]=\Q[d]$. De même, $\mathrm{Gr}^W_{-1} j_{!*}H'[d]=IC(H_{\Q}),$ ce qui achève de prouver l'exactitude du complexe ci-dessus.

\bigskip

Il est maintenant facile de conclure la preuve. La classe de l'extension
\begin{equation}\label{extperv}
0\ra IC(H_{\Q})\ra j_{!*}H'[d]\ra\Q[d]\ra 0
\end{equation}
donne, en localisant en tout point $b$ de $B$ un élément $\sigma$ de $IH^1(H_{\Q})_b$. D'autre part, son tiré en arrière par $j^*$ est l'extension venant de $\nu$. Il en résulte que l'image de $\sigma$ dans $H^1(H_{\Q})_b$ est $\sigma_b(\nu)$, ce qui conclut.\qed

\subsection{La conjecture de Hodge et la famille universelle des sections hyperplanes d'une variété}

Avant d'énoncer le résultat principal de cette section, nous aurons besoin d'une généralisation de \ref{ZP}. Soit $X$ une variété projective lisse complexe de dimension $2n$. Soit $\mathcal L$ un fibré en droite très ample sur $X$. Soit $\overline P$ l'espace projectif $|\mathcal L|$. On dispose de la variété d'incidence $\mathfrak X$ associée à la paire $(X, \mathcal L)$, dont les points sont les couples $(x, f)$ où $x$ est un point de $X$ et $f$ un point de $\overline P$ tel que $f(x)=0$. Soient $p$ et $\pi$ les projections de $\mathfrak X$ sur $X$ et $\pi$ respectivement. La flèche $p : \mathfrak X\ra X$ est un fibré projectif. 

Soit $X^{\vee}$ la variété duale de $X$, et soit $P$ le complémentaire de $X^{\vee}$ dans $\overline P$. Le morphisme $\pi$ est lisse au-dessus de $P$ par définition. Notons $H$ la variation de structures de Hodge de poids $-1$ sur $P$ associée à $R^{2n-1}\pi_*\Z(n)$, et $J\ra P$ le quotient de la fibration en jacobiennes intermédiaires $J(H)$ par la jacobienne intermédiaire $J^n(X)$ de $X$. Le procédé de \ref{ZP} permet, si $\alpha$ est une classe de Hodge dans $H^{2n}(X, \Z(n))$, d'associer à $\alpha$ une fonction normale $\nu_{\alpha}$ sur $J$. Les singularités des fonctions normales associées à la fibration constante $J(X)$ étant nulles comme on le vérifie sans difficulté, on consid\'erera les singularités $\sigma_{p}(\nu_{\alpha})$ de $\nu_{\alpha}$ comme des éléments de $IH^1(H_{\Q})_p$ grâce au théorème \ref{perv}.

\bigskip

Le but de cette section est d'esquisser la preuve du théorème suivant, d'après \cite{GG07, BFNP, dCM09}. Son sens est double. D'une part, bien que l'application d'Abel-Jacobi ne soit pas surjective en général, ce qui comme on l'a vu rend difficile l'application à la conjecture de Hodge de l'étude des fonctions normales associées aux pinceaux de Lefschetz, il montre que l'étude des fonctions normales au-dessus d'une base de grande dimension permet d'aborder la construction de cycles algébriques. D'autre part, il met en avant, par sa démonstration, l'importance dans le contexte de la conjecture de Hodge du théorème de décomposition de \cite{BBD}, et prolonge en ce sens le théorème \ref{perv}.

\begin{theo}\label{Hodge}
Les deux énoncés suivants sont équivalents.
\begin{enumerate}
\item La conjecture de Hodge vaut pour toutes les variétés projectives lisses complexes.
\item Si $X$ est une variété projective lisse complexe de dimension $2n$ munie d'un fibré en droites très ample $\mathcal L$, et si $\alpha$ est une classe de Hodge primitive dans $H^{2n}(X, \Z(n))$, on peut, quitte à remplacer $\mathcal L$ par $\mathcal L^k$ pour un entier $k$ assez grand, trouver un point $P$ de $\overline P$ tel que la singularité $\sigma_p(\nu_{\alpha})$ soit non nulle.
\end{enumerate}
\end{theo}

\noindent{\sc Preuve} ---
On va montrer que le deuxième énoncé implique le premier. Notons d'abord qu'un argument impliquant le théorème de Baire permet de réduire la conjecture de Hodge en degré $2n$ aux variétés de dimension $2n$. On raisonne par récurrence, et on suppose la conjecture de Hodge démontrée en tous degrés jusqu'à $2n-2$.

Commençons par appliquer le théorème de décomposition à $\pi$. On obtient une décomposition dans la catégorie dérivée
\begin{equation}\label{decomp}
R\pi_*\Q[2n]\simeq \bigoplus_i\,^pH^i(R\pi_*\Q[2n])[-i]
\end{equation}
du complexe $R\pi_*\Q$ en somme de ses faisceaux de cohomologie perverse. Soient $d$ la dimension de $P$ et $j$ l'inclusion de $P$ dans $\overline P$. Au-dessus de $P$, il s'agit de la décomposition donnée par la dégénérescence de la suite spectrale de Leray \cite{De68}, et l'on a $j^*\,^pH^i(R\pi_*\Q[2n])=R^{i+2n-d}(\pi_{|P})_*\Q[d]$.

Considérons le faisceau $^pH^i(R\pi_*\Q[2n])$. Il se décompose suivant le support strict. En particulier, il contient comme facteur direct l'extension intermédiaire de sa restriction à $P$, soit $IC(R^{i+2n-d}(\pi_{|P})_*\Q)$.

Ce qui précède permet d'exprimer les groupes 
$$H^{-i}(IC(R^{i+2n-d}(\pi_{|P})_*\Q))=IH^{-i+d}(R^{i+2n-d}(\pi_{|P})_*\Q)$$
comme facteurs directs de $H^{2n}(X, \Q)$. De même, si $p$ est un point de $\overline P$, la fibre de (\ref{decomp}) en $p$ permet d'exprimer les 
$$H^{-i}(IC(R^{i+2n-d}(\pi_{|P})_*\Q))_p=IH^{-i+d}(R^{i+2n-d}(\pi_{|P})_*\Q)_p$$ 
comme facteurs directs de $H^{2n}(\mathfrak X_p, \Q)$.

\bigskip

Bien entendu, ces expressions comme facteurs directs ne sont pas canoniques, car la décomposition (\ref{decomp}) ne l'est pas. Néanmoins, la flèche 
$$H^{2n}(X, \Q)\ra IH^0(R^{2n}(\pi_{|P})_*\Q)$$
est bien définie -- elle vient de la suite spectrale de foncteurs composés pour $\pi$ et le foncteur sections globales. Notons que le faisceau $R^{2n}(\pi_{|P})_*\Q$ est constant par le théorème de Lefschetz fort. La flèche ci-dessus correspond à la restriction d'une classe de cohomologie aux fibres lisses de $\pi$.

Si $\alpha$ est une classe de primitive comme dans le second énoncé du théorème, son image dans $IH^0(R^{2n}(\pi_{|P})_*\Q)$ est nulle par hypothèse. Pour la même raison, son image dans $IH^1(R^{2n-1}(\pi_{|P})_*\Q)$ est bien définie. On peut vérifier qu'il s'agit de la classe de l'extension (\ref{extperv}) dans la preuve du théorème \ref{perv}.

On peut appliquer le même raisonnement à la cohomologie de $\mathfrak X_p$. La restriction de la classe de $\alpha$ à $H^{2n}(\mathfrak X_p, \Q)$ a une image nulle dans le groupe $IH^0(R^{2n}(\pi_{|P})_*\Q)_p$, et définit par conséquent un élément canonique de $IH^1(R^{2n-1}(\pi_{|P})_*\Q)_p$. D'après ce qui précède, cet élément est la singularité $\sigma_p(\nu_{\alpha})$ de la fonction normale $\nu_{\alpha}$ en $p$.

\bigskip

Supposons la conjecture de Hodge fausse. Par le théorème de l'indice de Hodge et par le théorème de Lefschetz fort, on peut dans ce cas trouver une classe $\alpha$ qui soit d'intersection nulle avec tout cycle algébrique de dimension $n$ sur $X$.

Supposons maintenant que $\sigma_p(\nu_{\alpha})$ soit non nul pour un point $p$ de $\overline P$. Alors la restriction de $\alpha$ à $\mathfrak X_p$ est non nulle. Par hypothèse de récurrence, et quitte à travailler sur une résolution des singularités de $\mathfrak X_p$, cette restriction est une classe de cycle algébrique de codimension $n$ non nulle d'après la discussion ci-dessus. En particulier, on peut trouver sur une désingularisation de $\mathfrak X_p$ un cycle algébrique de dimension $n$ d'intersection non-nulle avec $\alpha_{|\mathfrak X_p}$. L'image de ce cycle dans $X$ a une intersection non nulle avec $\alpha$, d'où la contradiction cherchée.

\bigskip

Nous ne démontrons pas l'affirmation réciproque. Elle repose, outre les arguments précédents, sur le théorème de Thomas \cite{Th05} qui montre que la conjecture de Hodge implique que si $\alpha$ est une classe de Hodge non nulle dans $H^{2n}(X, \Q(n))$, on peut trouver une section hyperplane $\mathfrak X_p$ de $X$ telle que la restriction de $\alpha$ à $\mathfrak X_p$ est non nulle. Un théorème de Lefschetz faible pervers permet d'en déduire, si $\mathcal L$ est suffisamment ample, la non-nullité de la singularité de $\nu_{\alpha}$ en $p$.\qed

Il semble difficile de produire des singularités pour les fonctions normales. Un énoncé dans cette direction est le suivant. C'est une application simple du théorème \ref{algebrique} qui apparaît dans \cite{Sc10}.

\begin{prop}
Avec les notations ci-dessus, soit $\nu$ la fonction normale sur $J$ associée à une classe de Hodge primitive non nulle dans $H^{2n}(X, \Q(n))$. Si le lieu des zéros contient une composante de dimension strictement positive, $\nu$ est singulière en l'un des points où l'adhérence de cette composante rencontre la variété duale $X^{\vee}$.
\end{prop}

Il ne semble pas que ce dernier énoncé permette d'obtenir des résultats d'existence de cycles algébriques.

\section{Modèle de Néron pour les familles de jacobiennes intermédiaires}

Cette partie est consacrée à une description du modèle de Néron de Schnell introduit dans \cite{Sc09}, et à des éléments de preuve du théorème \ref{algebrique}.

\subsection{Modules de Hodge mixtes}

Comme entrevu dans la section précédente, l'étude des fonctions normales sur une base quelconque fait intervenir d'une part des faisceaux pervers et d'autre part des arguments de poids. La construction par Saito de la catégorie des modules de Hodge mixtes \cite{Sa88, Sa90} est le pendant de la théorie des faisceaux pervers $\ell$-adiques de \cite{BBD} au-dessus des corps finis et permet d'utiliser ce type d'arguments en théorie de Hodge. Nous rassemblons ici quelques éléments de cette construction.

Dans tout ce qui suit, nous travaillons sur une base $B$ analytique lisse, munie d'une compactification partielle lisse $j : B\hookrightarrow \overline B$. On note $d$ la dimension de $B$.

Un module de Hodge mixte est un objet qui généralise la notion de variation de structures de Hodge mixte. Sur $B$, un tel objet est constitué des données suivantes :
\begin{itemize}
\item un faisceau pervers K à coefficients rationnels, muni d'une filtration croissante $W_{\bullet}$ par des sous-faisceaux pervers,
\item un $\mathcal D_B$-module régulier holonome $\mathcal M$ muni d'une filtration croissante $W_{\bullet}$,
\item une bonne filtration croissante $F_{\bullet}$ de $\mathcal M$ par des sous-$\mathcal O_B$-faisceaux cohérents.
\end{itemize}

La filtration $W_{\bullet}$ est la filtration par le poids, et la filtration $F_{\bullet}$ est la filtration de Hodge.

Via la correspondance de Riemann-Hilbert \cite{Ka84, Me84}, on exige que le faisceau pervers associé à $\mathcal M$ soit $K_{\C}$ et que les filtrations par le poids se correspondent. Ces données sont soumises à des compatibilités supplémentaires définies par récurrence sur la dimension de $B$ à l'aide des faisceaux de cycles évanescents et de cycles proches. On dispose d'une notion de polarisation, et l'on ne consid\'erera que des modules de Hodge mixtes polarisables.

Un exemple de modules de Hodge mixte est donnée par une variation admissible $H$ de structures de Hodge mixtes sur $B$. Le faisceau pervers associé est $H_{\Q}[d]$, et le $\mathcal D_B$-module est $\mathcal H$ muni de sa connexion plate. La filtration de Hodge est dans ce cas donnée par $F_p=F^{-p}$. Pour cette raison, nous utiliserons désormais des filtrations de Hodge croissantes en faisant usage de la convention précédente.

Le point remarquable de la théorie de Saito est l'existence d'un formalisme des six opérations que nous utiliserons librement. Si $H$ est une variation de structures de Hodge mixtes sur $B$, admissible relativement à $\overline B$, on dispose ainsi de l'extension intermédiaire $M=j_{!*}H$. Le faisceau pervers sous-jacent à $M$ est obtenu par extension intermédiaire, et le $\mathcal D_B$-module par extension minimale. Le fait que $M$ soit un module de Hodge mixte est une des conséquences difficiles de la théorie de Saito \cite[3.10]{Sa90}.

Il est à noter que le $\mathcal D_B$-module sous-jacent à un module de Hodge mixte n'est en général pas cohérent comme $\mathcal O_B$-module. C'est le cas seulement s'il correspond à une variation de structures de Hodge mixte.

\subsection{Construction du modèle de Néron}

La construction de Schnell est l'exact analogue, dans le cadre des modules de Hodge mixtes, de la construction des jacobiennes intermédiaires en \ref{NF}. Soit $H$ une variation de structures de Hodge polarisable pure de poids $-1$, admissible relativement à $\overline B$. Le module de Hodge mixte correspondant est lui de poids $d-1$ -- remarquer que le faisceau pervers sous-jacent est $H_{\Q}$ placé en degré -d.

Soit $M$ le module de Hodge $j_{!*} M$. Il s'agit d'un module de Hodge polarisable sur $\overline B$, pur de poids $d-1$ lui aussi. Soit $\check M=\mathbb D_{\overline B}(M)(1-d)$ le dual de Verdier de $M$, normalisé pour avoir poids $d-1$. Sa restriction à $B$ est le module de Hodge mixte associé à la variation de structures de Hodge $H^{\vee}(1)$. On note $\mathcal M$ (resp. $\check{\mathcal M}$) le $\mathcal D_B$-module holonome régulier sous-jacent à $M$ (resp. $\check M$).

La filtration de Hodge fournit un faisceau $F_0 \check{\mathcal M}$ cohérent sur $\overline B$. Notons 
$$p : T(F_0\check{\mathcal M})\ra \overline B$$
le spectre de l'algèbre symétrique de $F_0 \check{\mathcal M}$. Il s'agit d'un espace analytique au-dessus de $\overline B$, dont le faisceau des sections est $(F_0 \check{\mathcal M})^{\vee}$. C'est en un certain sens l'espace total du faisceau cohérent $(F_0 \check{\mathcal M})^{\vee}$.

Au-dessus de $B$, on dispose d'un morphisme du système local $H_{\Z}$ dans $(F_0\check{\mathcal H})^{\vee}$. Soit maintenant 
$$p _{\Z} : T_{\Z}\ra \overline B$$
l'espace étalé du faisceau $j_*H_{\Z}$ sur $\overline B$. C'est un espace analytique muni d'un isomorphisme local vers $\overline B$ dont les sections s'identifient aux sections de $j_*H_{\Z}$.

Comme en \ref{NF}, l'inclusion de $H_{\Z}$ dans $(F_0\check{\mathcal H})^{\vee}$ au-dessus de $B$ correspond à un morphisme au-dessus de $B$
$$\epsilon_{|B} : (T_{\Z})_{|B} \ra T(F_0\check{\mathcal M})_{|B}.$$
La proposition \ref{epsilon} montre qu'il s'agit d'une immersion fermée d'espaces analytiques en groupes.

\begin{prop}\label{eps}
Le morphisme $\epsilon_{|B}$ se prolonge de manière canonique en un morphisme holomorphe de groupes au-dessus de $\overline B$
$$\epsilon : T_{\Z}\ra T(F_0\check{\mathcal M}).$$
\end{prop}

\noindent{\sc Preuve} --- 
Soit $K$ le faisceau pervers sous-jacent à $M$. Le faisceau $K$ est l'extension intermédiaire $j_{!*} H_{\Q}[d]$. En particulier, on a $\tau^{\leq -d} K\simeq H_{\Q}[d]$. D'autre part, par la correspondance de Riemann-Hilbert, le faisceau pervers $K_{\C}$ est quasi-isomorphe au complexe de de Rham
$$DR_{\overline B}(\mathcal M) = [\mathcal M\ra \Omega^1_{\overline B}\otimes \mathcal M \ra \ldots \ra  \Omega^d_{\overline B}\otimes \mathcal M][d].$$
Ces deux remarques montrent que l'on a un isomorphisme canonique de faisceaux 
sur~$\overline B$
$$j_*H_{\C}\ra \mathrm{Ker}(\mathcal M\ra \Omega^1_{\overline B}\otimes \mathcal M).$$
Autrement dit, le faisceau $j_*H_{\C}$ est le faisceau des sections plates de $\mathcal M$.

\bigskip

Un argument de dualité montre que l'on a un isomorphisme canonique de faisceaux sur $\overline B$
$$\mathrm{Ker}(\mathcal M\ra \Omega^1_{\overline B}\otimes \mathcal M)\simeq \mathcal Hom_{\mathcal D_{\overline B}}(\check{\mathcal M}, \mathcal O_{\overline B}).$$

Par bidualité et restriction, ces deux isomorphismes permettent de faire agir les sections de $j_*H_{\Z}$ sur les sections de $F_0\check{\mathcal M}$. On en déduit l'existence du morphisme $\epsilon$.\qed 

\begin{defi}
Le modèle de Néron de Schnell associé à la variation de structures de Hodge $H$ sur $B$ est le groupe topologique quotient au-dessus de $\overline B$
$$J^S(H)=T(F_0\check{\mathcal M})/\epsilon(T_{\Z}).$$
\end{defi}

Nous verrons plus bas que le morphisme $\epsilon$ est une immersion fermée. Cela permet de montrer que $J^S(H)$ est muni d'une structure canonique d'espace analytique séparé au-dessus de $\overline B$.

\bigskip

Un point tout à fait remarquable de cette définition est qu'elle est valable sans aucune hypothèse sur la géométrie du complémentaire de $\overline B$ dans $B$. Bien sûr, on peut supposer qu'il s'agit d'un diviseur car les variations de structures de Hodge s'étendent au-dessus des fermés de codimension au moins $2$.

Cela rend le modèle de Néron ci-dessus particulièrement adapté à l'étude de la famille universelle des sections hyperplanes d'une variété. Cette situation est abordée dans \cite{Sc12-2}. Décrivons brièvement un exemple concret. Je remercie Gérard Laumon de m'avoir communiqué ses notes sur la géométrie de la courbe plane universelle.

\bigskip

Soit $d$ un entier strictement positif, et soit $\mathcal C_0\ra \overline P$ la famille universelle des courbes planes de degré $d$. Une courbe de degré $d$ générique intersecte la droite à l'infini en $d$ points. Restreignons-nous aux courbes transverses à la droite à l'infini, et soit 
\mbox{$\pi : \mathcal C\ra U$} la courbe plane universelle obtenue en passant au revêtement de groupe~$\Sigma_d$ correspondant et en se restreignant à l'ouvert paramétrant les courbes réduites.

Soit $\mathcal J$ la composante neutre du champ de Picard relatif de $\mathcal C/U$. C'est un schéma en groupes commutatifs lisse de dimension relative $q=\frac{(d-1)(d-2)}{2}$. La restriction de $\mathcal J$ à l'ouvert de lissité $V$ de $\pi$ est bien entendu la famille des jacobiennes des fibres. Si $H$ est la variation de structures de Hodge de poids $-1$ associée à $R^1\pi_*\Z(1)$ au-dessus de l'ouvert $V$, alors la restriction de $\mathcal J$ à $V$ est la jacobienne intermédiaire $J(H)$. On peut montrer le résultat suivant.

\begin{prop}
L'espace analytique sous-jacent au schéma en groupes $\mathcal J\ra U$ est le modèle de Néron de Schnell associé à la famille de jacobiennes intermédiaires $J(H)$ sur $V\subset U$.
\end{prop}

\begin{rema}
Même dans cette situation simple, la géométrie des dégénérescences n'est pas facile à comprendre, et le complémentaire de $V$ dans $U$ est très singulier.
\end{rema}

\noindent{\sc Preuve} (Esquisse) ---
Appliquons le théorème de décomposition à $\pi$. La dimension de la base $U$ est $q+3d-1$. On obtient une décomposition
$$R\pi_*\Q[q+3d]=\bigoplus_i\,^pH^i(R\pi_*\Q[q+3d])[i].$$
Intéressons-nous comme dans \ref{Hodge} au faisceau pervers $^pH^0(R\pi_*\Q[q+3d])$. Sa restriction à l'ouvert de lissité $V$ est égale à
$$j^*\,^pH^0(R\pi_*\Q[q+3d])=R^1\pi_*\Q[q+3d-1].$$

Le théorème du support de Ngô \cite[7.2.1]{Ngo} -- que l'on pourrait remplacer dans ce cas par le théorème de Lefschetz pervers de \cite{BFNP} ou par \cite{Sc10} -- montre, après un argument de spécialisation, que le faisceau pervers $^pH^0(R\pi_*\Q[q+3d])$ admet $U$ comme support strict. Autrement dit, on a 
$$^pH^0(R\pi_*\Q[q+3d])=j_{!*}R^1\pi_*\Q[q+3d-1]=IC(R^1\pi_*\Q).$$

Le théorème de décomposition ci-dessus est encore valable dans la catégorie des modules de Hodge mixtes, et l'égalité ci-dessus vaut dans cette catégorie par pleine fidélité de la correspondance de Riemann-Hilbert. Cela montre que, dans la catégorie des modules de Hodge mixtes, on a un isomorphisme canonique
$$H^0(R\pi_*\Q[q+3d])=j_{!*} H$$
où $H$ est la variation de structure de Hodge pure de poids $-1$ associée à $R^1\pi_*\Q(1)$. L'égalité précédente permet de relier le modèle de Néron de Schnell, construit à partir du module de Hodge de droite, à la géométrie de $\pi$. Il n'est pas très difficile d'en déduire le résultat annoncé.\qed

\subsection{Séparation du modèle et extension des fonctions normales}

Nous décrivons maintenant, sans preuves, les principales propriétés du modèle de Néron construit ci-dessus. On garde les notations de la partie précédente. Tous les résultats sont tirés de \cite{Sc09}.

\begin{prop}\label{propre}
Le morphisme 
$$\epsilon : T_{\Z}\ra T(F_0\check{\mathcal M})_{|B}$$
est une immersion fermée propre.
\end{prop}

Une conséquence immédiate de ce qui précède est la propriété fondamentale suivante.

\begin{theo}\label{separe}
Le modèle de Néron $J^S(H)$ est un espace analytique -- en particulier séparé -- sur $\overline B$.
\end{theo}

Cette propriété de séparation est celle qui ne vaut pas pour le modèle de Zucker. C'est pour cette raison que l'utilisation des $\mathcal D$-modules est plus efficace ici que celle de l'extension de Deligne, qui choisit de prolonger le fibré vectoriel en perdant la structure différentielle -- puisque la connexion de Gauss-Manin acquiert des singularités.

Le fait que le $\mathcal D_{\overline B}$-module apparaissant dans la définition de $J^S(H)$ n'est pas cohérent comme $\mathcal O_{\overline B}$-module n'est pas gênant car on peut ne travailler qu'avec un sous-fibré de Hodge qui lui est de type fini.

\bigskip

Ce qui précède assurant que le modèle de Néron $J^S(H)$ est un objet analytique raisonnable, on peut analyser les propriétés de prolongement des fonctions normales.

\begin{prop}\label{ns}
Soit $\nu : B\ra J(H)$ une fonction normale admissible relativement à $\overline B$. La fonction normale $\nu$ se prolonge en une section $\overline B\ra J^S(H)$ si et seulement les singularités de $\nu$ sont toutes nulles.
\end{prop}

En fait, on peut montrer plus généralement que les sections globales horizontales de $J^S(H)$, pour une notion d'horizontalité adéquate, sont exactement les fonctions normales admissibles non-singulières $\nu : B\ra J(H)$.

Le résultat précédent résulte des bonnes propriétés de fonctorialité de $J^S(H)$. Ce qui suit est plus difficile.

\begin{theo}\label{closure}
Soit $\nu : B\ra J(H)$ une fonction normale admissible relativement à $\overline B$. L'adhérence topologique du graphe de $\nu$ dans $J^S(H)$ est un sous-ensemble analytique de $J^S(H)$.
\end{theo}

On en déduit immédiatement le théorème \ref{algebrique} -- dans le cas au moins des variations de structures de Hodge de poids $-1$.

\subsection{\'Eléments de preuve}

Nous donnons pour conclure quelques éléments de preuve des résultats annoncés dans la section précédente. La propriété la plus simple est celle de prolongement des fonctions normales admissibles sans singularités.

\noindent{\sc Preuve de la proposition \ref{ns}} --- 
Par \cite{Sa96}, la fonction normale admissible $\nu$ correspond à une extension de modules de Hodge mixtes
\begin{equation}\label{MHM}
0\ra M\ra N\ra \Q[d]\ra 0.
\end{equation}

D'autre part, la nullité des singularités de $\nu$ signifie précisément que l'on a une extension de faisceaux de groupes abéliens sur $\overline B$
\begin{equation}\label{local}
0\ra j_*H_{\Z}\ra  j_*H'_{\Z}\ra  \Z\ra 0
\end{equation}
obtenue en poussant en avant l'extension correspondant à $\nu$ sur $B$.

Le raisonnement de la proposition \ref{ext} s'applique : relevant localement la section triviale du système local $\Z$ dans la suite exacte (\ref{local}) et dualisant la suite exacte de $\mathcal D_{\overline B}$-modules sous-jacente à (\ref{MHM}), on produit effectivement par bidualité et restriction des sections locales de $(F_0\check{\mathcal M})^{\vee}$ bien définies à des sections de $j_*H_{\Z}$ près, soit une section globale de $J^S(H)$. Cela conclut.\qed

Les preuves de la proposition \ref{propre} et du théorème \ref{closure} sont plus difficiles. On peut remarquer tout d'abord que les deux résultats sont analogues. Soit $\nu$ une fonction normale admissible comme dans le théorème \ref{closure}. Elle induit une extension de systèmes locaux sur $B$
$$0\ra H_{\Z}\ra H'_{\Z}\ra \Z\ra 0.$$
Soit $T_{\nu}$ le sous-espace de l'espace étalé du faisceau $H'_{\Z}$ correspondant aux sections s'envoyant sur $1$ par la suite exacte précédente. On a un morphisme au-dessus de $B$ de $T_{\nu}$ dans $T(F_0\check{\mathcal M})$ obtenu comme plus haut. Le théorème \ref{closure} revient essentiellement à montrer que l'adhérence de l'image de $T_{\nu}$ dans $T(F_0\check{\mathcal M})$ est un sous-ensemble analytique.

Si maintenant la fonction normale $\nu$ est identiquement nulle, $T_{\nu}$ s'identifie à la restriction de $T_{\Z}$ à $B$ s'envoyant dans $T(F_0\check{\mathcal M})$ par $\epsilon$. \`A la propriété d'injectivité près, les énoncés de \ref{propre} et \ref{closure} se correspondent.

\bigskip

On considère maintenant la proposition \ref{propre} seulement. C'est la fonctorialité de $J^S(H)$ qui permet de se ramener à une situation géométrique accessible -- elle ne permet cependant pas par elle-même de montrer que le morphisme $\epsilon$ est une immersion fermée.

Montrons d'abord comment obtenir l'injectivité de $\epsilon$. Il suffit de raisonner fibre par fibre. Soit $i : \{b\}\hookrightarrow \overline B$ l'inclusion d'un point dans $\overline B$. Soit $H$ la structure de Hodge mixte $H^{-d}i^* M$. Il s'agit d'une structure de Hodge mixte dont tous les poids sont inférieurs à $-1$ car le poids de $M$ est $d-1$. On peut montrer, via la théorie des cycles évanescents, qu'elle est munie d'une structure entière via l'inclusion de $T_{\Z, b}$ et que l'on a une surjection 
$$T(F_0\check{\mathcal M})_b\ra H_{\C}/F_0H_{\C}$$
compatible aux inclusions de $T_{\Z, b}$. Comme les poids de $H$ sont inférieurs à $-1$, $T_{\Z, b}$ s'identifie à un sous-groupe discret de $H_{\C}/F_0H_{\C}$. Cela montre que la restriction de $\epsilon$ à $T_{\Z, b}$ est injective.

\bigskip

Soit $f : \overline C\ra \overline B$ une application holomorphe. On suppose que l'image réciproque $C$ de $B$ est dense dans $\overline C$. Soit $d_C$ la dimension de $C$. Soit $M^f$ le module de Hodge mixte 
$$M^f=H^{d-d_C} f^! M(d-d_C).$$
Cette formule traduit la construction suivante : le module de Hodge mixte $M^f$ est l'extension intermédiaire de $C$ à $\overline C$ du module de Hodge mixte associé au tiré en arrière à $C$ de la variation de structures de Hodge $H$ sur $B$.
La fonctorialité pour les modules de Hodge mixtes donne un morphisme canonique
$$F_0\check{\mathcal M}^f\ra f^*F_0\check{\mathcal M}.$$
C'est un isomorphisme au-dessus de $C$. Suivant ces compatibilités, on obtient la fonctorialité des modèles de Néron sous la forme d'un morphisme holomorphe canonique
$$J^S(H)\times_{\overline B}\overline C \ra J^S((f_{|C})^*H).$$

Par un raisonnement semblable à celui que nous avons utilisé pour montrer l'injectivité, la fonctorialité ci-dessus permet de ramener la preuve de la proposition \ref{propre} -- et celle du théorème \ref{closure} -- pour $H$ à celle pour $f^*H$ si $f$ est propre.

Après résolution des singularités du complémentaire de $B$ dans $\overline B$, passage à un revêtement fini et restriction à un ouvert, cela permet de se ramener au cas où $B$ est un produit $(\Delta^*)^d)$ de disques épointés inclus dans $\overline B=\Delta ^d$.

\bigskip

Dans cette situation, l'extension intermédiaire $M$ admet une description explicite en termes de l'extension de Deligne $\overline{\mathcal H}$ de $\mathcal H$. L'extension de Deligne est un sous-faisceau de $j_*\mathcal H$. Ce dernier est muni d'une structure de $\mathcal D_{\overline B}$-module venant de la connexion de Gauss-Manin sur $\mathcal H$. Le $\mathcal D_{\overline B}$-module $\mathcal M$ est alors le $\mathcal D_{\overline B}$-module engendré par $\overline{\mathcal H}$. Autrement dit, les sections de $\mathcal M$ sont les différentielles itérées des sections de $\overline{\mathcal H}$. La filtration de Hodge de $\mathcal M$ est la convolution de la filtration venant de l'extension des fibrés de Hodge de $\mathcal H$ à $\overline{\mathcal H}$ et de la filtration de $\mathcal D_{\overline B}$ par l'ordre des formes différentielles.

C'est cette description qui explique pourquoi le modèle de Néron de Schnell est séparé tandis que l'extension de Zucker ne l'est pas. En effet, montrer que le morphisme $\epsilon$ est une immersion fermée revient à montrer que les sections de $F_0\check{\mathcal M}$, via la polarisation, séparent uniformément les sections de $j_*H_{\Z}$ au voisinage des points du complémentaire de $B$ dans $\overline B$. D'après ce qui précède, il y a plus de sections de $F_0\check{\mathcal M}$ que de sections de $\overline{\mathcal H}$. 

C'est ce point qui permet à Schnell de démontrer la proposition \ref{propre} et le théorème \ref{closure}. La mise en \oe uvre technique de cette idée repose sur des estimées délicates venant de la théorie de Schmid, et notamment d'un théorème de l'orbite SL$_2$ prouvé dans \cite{KNU08}. Il s'agit ici de borner les normes de certaines classes de cohomologie réelles, et ce sont ces bornes qui expriment le caractère propre de $\epsilon$.

\end{document}